# A REDUCTION METHODOLOGY USING FREE-FREE COMPONENT EIGENMODES AND ARNOLDI ENRICHMENT

**Hadrien Tournaire**
IRT SystemX
Palaiseau, France

**Franck Renaud**
LISMMA
Saint-Ouen, France

**Jean-Luc Dion**
LISMMA
Saint-Ouen, **France**

## ABSTRACT

In order to perform faster simulations, the model reduction is nowadays used in industrial contexts to solve large and complex problems. However, the efficiency of such an approach is sometimes cut by the interface size of the reduced model and its reusability.

In this article, we focus on the development of a reduction methodology for the build of modal analysis oriented and updatable reduced order model whose size is not linked to their contacting interface. In order to allow latter model readjusting, we impose the use of eigenmodes in the reduction basis. Eventually, the method introduced is coupled to an Arnoldi based enrichment algorithm in order to improve the accuracy of the reduced model produced.

In the last section the proposed methodology is discussed and compared to the Craig and Bampton reduction method. During this comparison we observed that even when not enriched, our work enables us to recover the Craig and Bampton accuracy with partially updatable and smaller reduced order model.

## INTRODUCTION

The iterative design of a system involves several validation simulations that may slow its conception process. Moreover, in such a context, the required accuracy and detail level is expected to provide a model whose significant size leads to time consuming simulations.

Nowadays substructuring and reduction methods like the superelements or the Craig and Bampton method as instances are industrially used to solve large and complex structural dynamic problems and to cut simulation time.

In a design context, the Mac Neal and Rubin reduction methods are particularly interesting since they use substructures eigenmodes. Indeed, the use of free-free component modes enables to evaluate the contribution of each part into the global dynamics of a system. This may contain interesting information for the design of a system in terms of dynamic behavior. Furthermore, using eigenmodes enables to build updatable models that do not have to be entirely rebuilt when the relative position of a component is changed. Eventually, the use of these modes also allows to perform model readjusting.

Despite their numerous qualities, the Mac Neal and Rubin methods such as Craig and Bampton are facing limitations since they lead to reduced models whose size depends on the substructures contacting interfaces. As a results, this can lower the interest of a reduction approach when working on heavy models.

Several researches have been led to overcome these weaknesses. As an example, methodologies to improve the reusability of the Craig end Bampton basis by enrichment were performed [1], [2]. Rixen proposed a dual Craig and Bampon method using free-free eigenmodes [3] and an additional contribution that aimed at reducing the number of interface mode [5].

We suggest here an alternative methodology providing reduced order model whose size is not linked to the contacting interface. The research of reusability for the reduced order model and the use of free-free eigenmodes in the reduction basis are conjointly performed. In a last stage, an Arnoldi based enrichment algorithm for a selective improvement of the model solutions is used.

## NOMENCLATURE

$\Sigma_k$ : Component k
$\Gamma_k$ : Contacting surface of the k component

$j_k$ : Junction DoF set of the k component
$i_k$ : Interior DoF set of the k component k

$M_k$ : Mass matrix of the k component
$K_k$ : Stiffness matrix of the k component
$Z_k$ : Dynamic stiffness of the k component

φ   : Eigenmodes basis (truncated or not)
φ_c : Basis of coupling vectors
T   : Reduction basis

SVD : Singular Value Decomposition

$\varepsilon_{tol}$ : Tolerance threshold
r   : Root of the Rayleigh coefficient

**Quality criterion of a reduced order model**

The quality estimation of a reduced order model can be done thanks to several tools, some of them are summarized in [8]. The selection of the quality criterion should be done in regard to the aim of the simulation.

Here the reduction is performed to faster the modal analysis of a structure, so that the Modal Assurance Criterion (MAC) is favored. Indeed, this indicator is geometric and enables to quantify the collinearity of two vectors. When dealing with vectors families $U = [u_1 \ ... \ u_n]$ and $V = [v_1 \ ... \ v_m]$ the terms of the associated MAC matrix of dimension $n \times m$ are defined as:

$$MAC_{i,j} = mac(u_i, v_j) = \frac{|u_i^T * v_j|^2}{u_i^T * \overline{u}_i * v_j^T * \overline{v}_j}$$

Where $\overline{u}$ is the conjugate of a vector u. Thus the quality of each solution of the reduced model is estimated by a MAC comparison with corresponding solutions of the non-reduced model.

**Reduction methodology**

Let us consider a system $\Sigma$ made of two components $\Sigma_1$ and $\Sigma_2$ tied to each other on their contacting interfaces $\Gamma_1$ and $\Gamma_2$ of the two components. For the sake of simplicity the two meshes of the contacting interfaces are taken as coherent and the interfaces DoF $j_1$ and $j_2$ are sorted so that they exactly match in both FE characteristic matrix. The interior DoF of the components are noted $i_1$ and $i_2$.

In order to get coherent displacement fields on the contacting interface the first step of our work consists in primarily assemble both components so that the deformation continuity is imposed in a strong manner. Indeed, dual assemblies rely on the force continuity that may cause displacement jumps on the contacting interfaces. The dynamic stiffness of the assembled system $Z(\omega)$ can be split into an elastic dynamic stiffness $Z_E(\omega)$ that is block diagonal and an interaction stiffness matrix $Z_I(\omega)$ that is mostly sparse so that:

$$Z(\omega) = Z_E + Z_I$$

$$Z(\omega) = \begin{bmatrix} Z_{1,ii} & Z_{1,ij} & 0 & 0 \\ Z_{1,ji} & Z_{1,jj} & 0 & 0 \\ 0 & 0 & Z_{2,ii} & Z_{2,ij} \\ 0 & 0 & Z_{2,ji} & Z_{2,jj} \end{bmatrix} + \begin{bmatrix} 0 & 0 & 0 & 0 \\ 0 & Z_{2,jj} & Z_{2,ji} & 0 \\ 0 & 0 & 0 & 0 \\ Z_{1,ji} & 0 & 0 & Z_{1,jj} \end{bmatrix}$$

The reduction of such a model with a basis $\phi_f$ only made of free-free eigenmodes $\phi_1$ and $\phi_2$ lead to reduced model exclusively working for very low frequencies [6]. Thus, additional reduction vectors $\phi_c$ called coupling vectors have to be taken into account to enable recovering the system behavior on a more spread frequency bandwidth. From now, the hunted reduction basis T have the following structure:

$$T = [\phi_f \ \phi_c]$$

In order to increase the numerical conditioning of the reduced order model the coupling deformations $\phi_c$ are orthogonalized to the free-free eigenmodes $\phi_f$ using a Gram-Schmidt algorithm.

**Interface reduction using the eigenmodes interface deformations**

In the previous section the use free-free eigenmodes $\phi_f$ to reduce the model was not satisfying, indeed this involves no global deformation of $\Sigma_1$ and $\Sigma_2$ simultaneously. Moreover, this projection of the coupled components on an uncoupled deformation basis may lead to interface discontinuities. In order to improve the displacement field accuracy on the interface a first method consist in imposing to each component the interface displacement field that are contained in the eigenmodes of the other component. The deformation of a k-component $\Sigma_k$ whose contacting interface $\Gamma_k$ is submitted to the deformations of a p-component $\Sigma_p$ are written $\Theta_{k/p}$:

$$T = \begin{bmatrix} \phi_{1,i} & 0 & \Theta_{1/2,i} & 0 \\ \phi_{1,j} & 0 & \Theta_{1/2,j} & 0 \\ 0 & \phi_{2,i} & 0 & \Theta_{2/1,i} \\ 0 & \phi_{2,j} & 0 & \Theta_{2/1,j} \end{bmatrix}$$

Where the deformations $\Theta_{k/p,j}$ are in fact equal to the interface deformation found in the eigenmodes of the p-component $\phi_{p,j}$. The deformations on the interior DoF of the components are computed using:

$$\begin{cases} \Theta_{1/2,i} = -Z_{1,ii}^{-1} * Z_{1,ij} * \phi_{2,j} \\ \Theta_{1/2,i} = -Z_{2,ii}^{-1} * Z_{2,ij} * \phi_{1,j} \end{cases}$$

In order to improve validity of the reduced order model it is interesting to use several coupling deformations computed for different given circular frequencies $\omega_1, \omega_2, ..., \omega_n$. In spite of this, the size of the basis T may fast increase and its conditioning may collapse, especially if the circular frequencies $\omega_1, \omega_2, ..., \omega_n$ are closed.

The coupling vectors presented below enable to get interesting results in terms of MAC but the associated reduced model may be consequent in terms of size, depending on the number of truncated modes in the basis $\phi_1$ and $\phi_2$. In order to increase the compactness of the basis T we use the Rayleigh coefficient to filter the coupling deformations that are remote from the studied frequency bandwidth.

The Rayleigh coefficient root's r of a coupling deformation $\Theta_{k/p}$ is computed using the characteristic matrix associated to the k-component, so that:

$$r = \sqrt{\frac{\Theta_{k/p}^T \ K_k \ \Theta_{k/p}}{\Theta_{k/p}^T \ M_k \ \Theta_{k/p}}}$$

Eventually, the coupling deformations are selected in regard to their Rayleigh coefficient. Indeed, the more a Rayleigh coefficient closed to the studied waveband is, the more the associated coupling vector is interesting for the model reduction.

**Interface reduction using a singular value decomposition of the eigenmodes interface displacements (proposed methodology)**

Another approach for the interface reduction consists in building a common interface displacement space $T_j$ for both components and simultaneously imposing these deformations to both components $\Theta_{T_j}$. We suggest here the use of the singular value decomposition (SVD) for the build of the interface displacement basis:

$$T_j = SVD(\phi_{1,j}, \phi_{2,j})$$

Like in the previous method, the coupling deformations are computed for different circular frequencies using the analog formula:

$$\begin{cases} \Theta_{T_j,1} = -Z_{1,ii}^{-1} * Z_{1,ij} * T_j \\ \Theta_{T_j,2} = -Z_{2,ii}^{-1} * Z_{2,ij} * T_j \end{cases}$$

Once again, it is interesting to compute these deformations for several pulsations. Here the selection of the coupling deformations is performed according to the singular value associated to the interface deformations computed with the SVD.

Indeed the higher singular values are associated to the deformations that best engender the interface deformations observed in the eigenmodes.

**Estimators for the quality of reduce model solutions**

Once built, the quality of a reduced order model has to be estimated. In our context the use of the MAC is not possible since we try to pass by the heavy computation of the high fidelity model (HFM) solutions.

Several indicators were developed for an a priori accuracy evaluation of a model solution, some of them are presented in the following works [7], [8]. In our case we decided to use the following force based indicators:

$$\varepsilon_i = \frac{\|R_f(\widetilde{\phi}_i)\|_K^2}{\|K\widetilde{\phi}_i\|_K^2} \text{ using the scalar product } \|g\|_K = g^T K g$$

Where $R_f(\widetilde{\phi}_i)$ is the residual force associated to the $i^{th}$ solution couple $(\widetilde{\omega}_i^2, \widetilde{\phi}_i)$ of the considered reduced model. Indeed, contrary to the displacement based indicators, a time-expensive inversion of the stiffness matrix K is avoided here.

Each solution $\widetilde{\phi}_i$ whose indicators reach a given tolerance $\varepsilon_{tol}$ is considered as poor and then enriched by an Arnoldi based algorithm that leads to construction of the new basis T. Once enriched, the solutions of the newly reduced model are computed and their quality indicator are evaluated once again.

In practice, each poor solution is enriched using n Arnoldi vectors that may fast lead to consequent reduction basis. Thus, the evaluations of the model solutions are expected to slow since the size of the reduce models tend to increase. To avoid this we use a restart step that consists in using the model solution in addition to the Arnoldi vectors to reduce the model Figure. 1.

Eventually, the coupling vectors computed in the loop are orthogonalized to the free-free eigenmodes and added to them that gives the reduction basis T.

**Update of the reduction basis**

Unlike the Craig and Bampton method, the modification of the relative position of the contacting surface has a limited impact on the reduction basis. Indeed, it only affects the coupling deformations. This explains why it is interesting to also build a database with containing both coupling and free-free deformations $\phi_l$, $\Theta_{k/p,j}$ and/or $\Theta_{T_j}$ of the component in order to faster the rebuild of the initial reduction basis T and then enrich it.

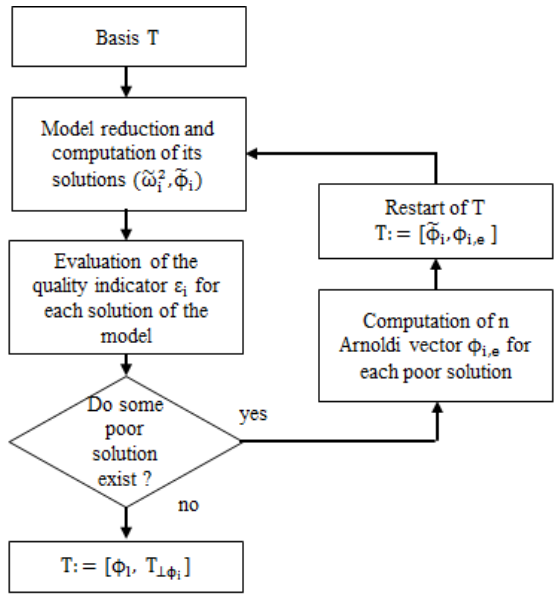

Figure 1.: Enrichment algorithm implementation

It has however to be noticed that in the case of a structural modification of a component, the proposed methodology do not offer any advantage compared to the Craig and Bampton method.

**Application of the reduction methodology**

In order to support our purpose we apply the methodology previously described to a study case made of two simple parts. The parts of the model are linear and considers small deformations, their properties can be associated to those of a steel. The two parts are embedded and studied on a bandwidth $\Delta f$ between 0 and 3000 Hertz that contains 41 eigenmodes.

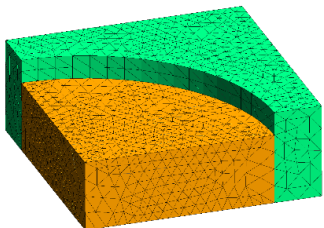

Figure 2.: Case study used for the test of the proposed methodology. the contacting interface of the component is curved and own 109 nodes (321 DoF)

The two interface reduction methods are compared themselves and with the Craig and Bampton method since it is a common method often considered as a reference in substructuring and that is already implemented in industrial codes like NASTRAN or CODE_ASTER.

Application of the Craig and Bampton method

Here, the tied interface eigenmodes used are taken in the frequency range between [0; 4000] Hertz (54 modes) while the interface deformations are computed for a null pulsation (321 vectors).

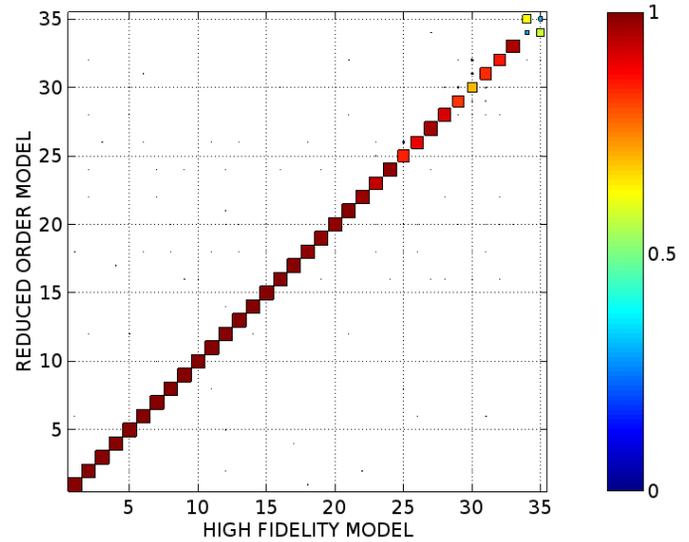

Figure 4.a.: MAC comparison of the flexible model solutions of the HFM and Craig and Bampton reduced model. Number of DoF: 375, mac average: 93.98%. Reduction basis conditioning: 133.65.

Application of the interface reduction methods

For the interface reduction methods, the free-free eigenmodes considered are taken below 3000 Hertz (47 modes). In order to get a rich model we computed the coupling deformations at the start, middle and end of the bandwidth (0, 1500 and 3000 Hertz).

The coupling deformations obtained using the SVD are selected by associated singular value. On Figure 3, the plot of the singular values of the interface deformation led us to only consider the singular values equal or higher than 0.2, this eventually lead to 126 coupling deformations.

In order to get reduction basis of the same size for the two methods, only the coupling vectors associated to the 126[th] smallest Rayleigh coefficients are taken into account for the firstly introduced interface reduction method.

As is it visible on Figure 4.b, the first interface reduction method provides rather good results for most of MAC but its conditioning is extremely poor that avoid to recover some of the high fidelity modes. This explains why its flexible MAC average collapses to 86.63% only. Indeed, the interface deformations observed in the eigenmodes of both components are rather similar. Moreover the pulsations taken for the computation of the coupling vectors are not remote enough to insure a good orthogonality between these vectors.

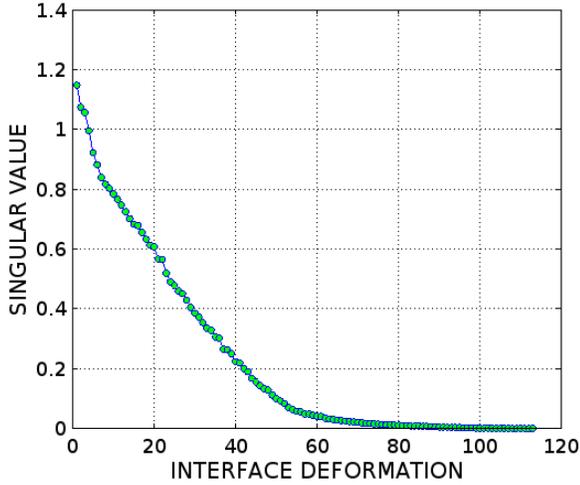

Figure 3.: Decrease of the singular values obtained by the SVD of the interface deformation observed in the eigenmodes.

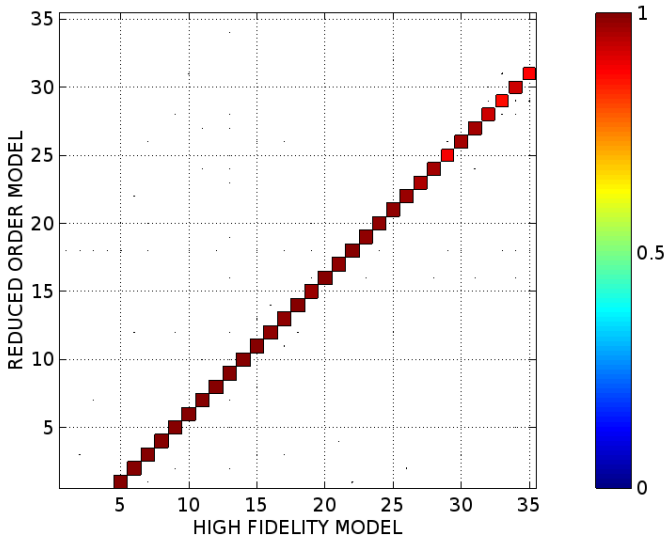

Figure 4.b. MAC comparison of the flexible model solutions of the HFM and reduced model using $\theta_{k/p,j}$ as coupling deformations. Number of DoF: 173, MAC average: 86.63%. Reduction basis conditioning: 1.29e+15.

Nevertheless, we can observe on Figure 4.c that the method using the singular decomposition of the interface deformations provides very good results in terms of flexible MAC average and compactness. Indeed, compared to the Craig and Bampton method the compactness of the reduction basis is increase by 2 while the model accuracy is improved by 2.44% of MAC average. It can however be noticed that the conditioning of the reduction basis is higher (9069.67) than the one of the Craig and Bampton method (133.65). This loss is mostly brought by the computation of the coupling deformations for too closed frequencies (0, 1500 and 3000 Hz).

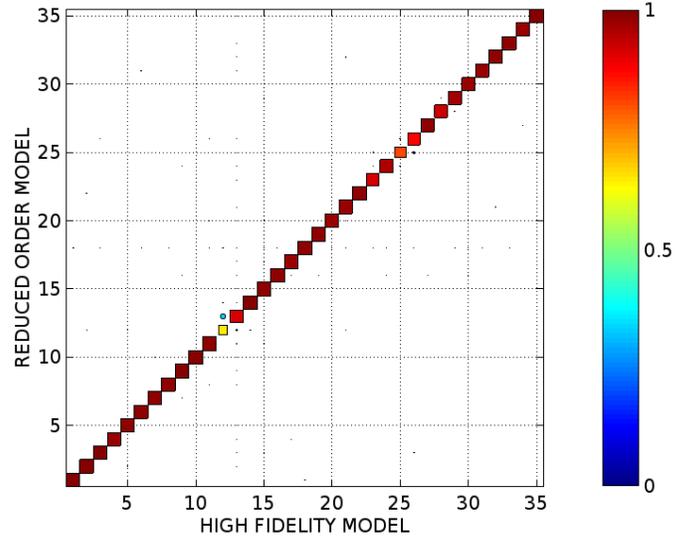

Figure 4.c. MAC comparison of the flexible model solutions of the HFM and reduced model using $\theta_{t_j}$ as coupling deformations. Number of DoF: 173, mac average: 96.42%. Reduction basis conditioning: 9069.67.

### Conclusions

The proposed method enable us to get accurate reduction basis for the Arnoldi enrichment. Thus, we expect this enrichment to be quick. The reduced model we build is more compact than those obtained with Craig and Bampton that was one of our first goal.

Although the reduced order model provided by the proposed reduction model is efficient for modal analysis it has to be noticed that its spatial convergence is not guaranteed for forced response simulations, unlike the Craig and Bampton method.

Eventually, the topology of the reduction basis we created is favorable to update if the position of a component is changed.

### ACKNOWLEDGMENTS

This research work has been carried out under the leadership of the Technological Research Institute SystemX, and therefore granted with public funds within the scope of the French Program "Investissements d'Avenir". In addition, we acknowledge the SNECMA and the LISMMA for their implication and cooperation on the current work.